\documentclass[reqno,12pt]{amsart}
\usepackage{amsmath,amsrefs,amssymb}

\numberwithin{equation}{section}

\DeclareMathOperator{\divergence}{div}
\DeclareMathOperator{\loc}{loc}

\DeclareMathOperator{\sgn}{sgn}
\DeclareMathOperator{\bigO}{O}
\DeclareMathOperator{\smallo}{o}
\newcommand{\R}{\mathbb{R}}
\renewcommand{\S}{\mathbb{S}}
\newcommand{\N}{\mathbb{N}}
\newcommand{\<}{\left<}
\renewcommand{\>}{\right>}
\renewcommand{\[}{\left[}

\renewcommand{\(}{\left(}
\renewcommand{\)}{\right)}

\newtheorem{theorem}{Theorem}[section]
\newtheorem{definition}{Definition}[section]
\newtheorem{proposition}{Proposition}[section]
\newtheorem{step}{Step}[section]
\newtheorem{remark}{Remark}[section]

\begin{document}

\title[Weighted critical $p$--Laplace equations]{Sharp pointwise estimates for weighted critical $p$--Laplace equations}

\author{Shaya Shakerian}

\address{Shaya Shakerian, University of California Santa Barbara, Department of Mathematics, CA 93106, USA}
\email{shaya@math.ucsb.edu}

\author{J\'er\^ome V\'etois}

\address{J\'er\^ome V\'etois, McGill University, Department of Mathematics and Statistics, 805 Sherbrooke Street West, Montreal, Quebec H3A 0B9, Canada}
\email{jerome.vetois@mcgill.ca}

\thanks{To appear in {\it Nonlinear Analysis : Theory, Methods \& Applications}.}

\thanks{The second author was supported by the Discovery Grant RGPIN-2016-04195 from the Natural
Sciences and Engineering Research Council of Canada.}

\date{March 5, 2020}


\begin{abstract}
We investigate the asymptotic behavior of solutions to a class of weighted quasilinear elliptic equations which arise from the Euler--Lagrange equation associated with the Caffarelli--Kohn--Nirenberg inequality. We obtain sharp pointwise estimates which extend and improve previous results obtained in the unweighted case. In particular, we show that we can refine the asymptotic expansion at infinity by using a Kelvin-type transformation, which reduces the problem to another elliptic-type problem near the origin. The application of this transformation is straightforward in the linear case but more delicate in the quasilinear case. In particular, it is necessary in this case to establish some preliminary estimates before being able to apply the transformation.
\end{abstract}

\maketitle

\section{Introduction and main results}

In this paper, we are interested in the elliptic problem
\begin{equation}\label{Main Problem}
\left\{\begin{aligned}& - \divergence \(\left|x\right|^{-ap}\left|\nabla u\right|^{p-2} \nabla u\) =f\(x,u\)&&\text{in }\R^n\\
&u\in D^{1,p}\(\R^n,\left|x\right|^{-ap}\),
\end{aligned}\right.
\end{equation}
where $f:\R^n\times\R\to\R$ is a Caratheodory function satisfying \begin{equation}\label{Th1Eq1}
\left|f\(x,s\)\right|\le\Lambda\left|x\right|^{-bq}\left|s\right|^{q-1}\quad\text{for all }s\in\R\text{ and a.e. }x\in\R^n
\end{equation}
for some constant $\Lambda>0$, the numbers $p$, $a$, $b$ and $q$ are such that
\begin{equation}\label{K-C-N ineq:conditions}
p>1,\quad a<\frac{n-p}{p}\,,\quad a\le b<a+1,\quad q=\frac{np}{n-p\(1+a-b\)}
\end{equation}
and 
$D^{1,p}(\R^n, \left|x\right|^{-ap})$ is the completion of $C^\infty_c\(\R^n\)$ with respect to the norm 
$$\left\|u\right\|_{D^{1,p}\(\R^n,\left|x\right|^{-ap}\)}:=\(\int_{\R^n}\left|x\right|^{-ap}\left|\nabla u\right|^pdx\)^{\frac{1}{p}}.$$
When $f\(x,u\)=\Lambda\left|x\right|^{-bq}\left|u\right|^{q-2}u$, (\ref{Main Problem}) is the Euler--Lagrange equation associated with the Caffarelli--Kohn--Nirenberg inequality~\cite{CafKohnNir}, which states, under conditions (\ref{K-C-N ineq:conditions}), that there exists a positive constant $C=C\(n,a,b,p\)>0$ such that
\begin{equation}\label{K-C-N ineq}
\(\int_{\R^n}\left|x\right|^{-bq}\left|u\right|^q dx \)^{\frac1q} \le C \( \int_{\R^n}\left|x\right|^{-ap}\left|\nabla u\right|^p dx \)^{\frac1p}
\end{equation}
for all $u\in C^\infty_c\(\R^n\)$.

\smallskip
Our main result is the following:

\begin{theorem}\label{Th1}
Let $n$, $p$, $a$, $b$ and $q$ be such that \eqref{K-C-N ineq:conditions} holds true, $f:\R^n\times\R\to\R$ be a Caratheodory function satisfying \eqref{Th1Eq1} and $u$ be a solution of \eqref{Main Problem}. Then there exists a constant $C_0>0$ such that 
\begin{equation}\label{Th1Eq3}
\left|x\right|^\mu\left|u\(x\)\right|+\left|x\right|^{\mu+1}\left|\nabla u\(x\)\right|\le C_0\quad\forall x\in \R^n\backslash B\(0,1\),
\end{equation} 
where $\mu:=\(n-p\(1+a\)\)/\(p-1\)$. If moreover $u>0$ and $f\(x,u\)\ge0$ in $\R^n$, then there exist constants $\alpha,\delta,C_1>0$ such that
\begin{equation}\label{Th1Eq4}
\left|\left|x\right|^\mu u\(x\)-\alpha\right|\le C_1\left|x\right|^{-\delta}\quad\forall x\in \R^n\backslash B\(0,1\)
\end{equation}
and 
\begin{equation}\label{Th1Eq5}
\left|x\right|^{\mu+1}\nabla u\(x\)+\alpha\mu\left|x\right|^{-1}x\to0\quad\text{as }\left|x\right|\to\infty.
\end{equation}
\end{theorem}

This theorem extends previous results obtained by Sciunzi~\cite{Sci} and V\'etois~\cite{Vet} in the case where $a=b=0$. In fact, for positive solutions, \eqref{Th1Eq4} and \eqref{Th1Eq5} improve the estimates obtained in~\cites{Sci,Vet}, where it was proven that $1/C\le\left|x\right|^\mu u\(x\),\left|x\right|^{\mu+1}\left|\nabla u\(x\)\right|\le C$ for some constant $C>0$ independent of $x\in\R^n\backslash B\(0,1\)$.

\smallskip
The proof of \eqref{Th1Eq4} relies in particular on the use of a Kelvin-type transformation of the form $u_\ast\(y\):=\left|x\right|^{\mu}u\(x\)$, where $x:=\left|y\right|^{-2}y$. This transformation is well-known in the case $p=2$ but, as far as the authors know, it has never been used in the case where $p\ne2$. In contrast with the case where $p=2$, the equation is not invariant under this transformation when $p\ne2$. In this paper, we show that we can still apply this transformation in the case where $p\ne2$, however it is necessary in this case to establish some preliminary estimates (see \eqref{Firstorder} below) and then the transformation can be used to improve these estimates and obtain for instance the H\"older-type estimate in \eqref{Th1Eq4}. We believe that even sharper results could be achieved by this method if new regularity results were obtained for the weighted elliptic-type problem that results from this transformation.

\smallskip
We point out that in contrast with the case where $a=b=0$, the solutions of \eqref{Main Problem} are not always radial in the presence of weights. Indeed, Horiuchi~\cite{Hor} obtained the existence of extremals when $a<0$ and $a<b<a+1$ (see also Catrina and Wang~\cite{CatWang}). However, these extremals turn out to be non-radial when $a<b<h\(a\)$ for some function $h$ such that $a<h\(a\)<a+1$ (see Catrina and Wang~\cite{CatWang} and Felli and Schneider~\cite{FelSch1} for $p=2$ and Byeon and Wang~\cite{ByeWang}, Caldiroli and Musina~\cite{CalMus} and Smets and Willem~\cite{SmeWil} for $p\ne2$). Caldiroli and Musina also observed in~\cite{CalMus}*{Section~4} that the Kelvin-type transformation $u_\ast\(y\):=u\(x\)$, where $x:=\left|y\right|^{-2}y$, transforms the problem into a similar weighted elliptic-type problem, where the only difference with the original problem lies in the exponents of the weights. We note that this transformation is however different from the one we use in this paper since we multiply here by the weight $\left|x\right|^{\mu}$ in order to remove the singularity at 0 in the transformed function. 

\smallskip
When $p=2$ and $0\le a\le b<a+1$, Chou and Chu~\cite{ChouChu} obtained that the positive solutions of \eqref{Main Problem} are always radial, thus extending the well-known result obtained by Caffarelli, Gidas and Spruck~\cite{CafGidSpr} for the classical Laplace operator i.e. when $a=b=0$. The optimal symmetry result in the case where $a<0$ and $p=2$ has recently been established by Dolbeault, Esteban and Loss~\cite{DolEstLoss1}. We also refer on this topic to the recent survey article~\cite{DolEstLoss2} by the same authors.

\smallskip
In the case where $a=b=0$ and $p\ne2$, the pointwise estimates obtained by Sciunzi~\cite{Sci} and V\'etois~\cite{Vet} have found application for establishing symmetry results for the solutions of \eqref{Main Problem}; see Damascelli and Ramaswamy~\cite{DamRam} and V\'etois~\cite{Vet} for $p<2$ and Sciunzi~\cite{Sci} for $p>2$. These results have recently been extended by Oliva, Sciunzi and Vaira~\cite{OliSciVai} to a class of $p$-Laplace equations with Hardy potential, using pointwise estimates established by Xiang~\cites{Xia1,Xia2} in this case. These results have also been extended, under a different method, still relying on pointwise estimates, by Ciraolo, Figalli and Roncoroni~\cite{CirFigRon} to a class of $p$-Laplace-type equations in an anisotropic setting. We also refer to the work by Esposito~\cite{Esp} which treats the limit case where $p=n$ and the nonlinearity is exponential, and where, again, pointwise estimates play a crucial role. 

\smallskip
The interest of this paper is therefore twofold. On the one hand, we believe that like in the aforementioned papers, our results will be useful to establish new symmetry results. On the other hand, since, as mentioned above, there exist situations where non-radial solutions exist, we are interested in developing a method to obtain sharp pointwise estimates for the solutions in this case.

\smallskip
We prove Theorem~\ref{Th1} through several steps. In Section~\ref{SecRegularity}, we obtain global boundedness results in $L^\infty\(\R^n\)$ and weak Lebesgue spaces. In Section~\ref{SecApriori1}, we establish \eqref{Th1Eq3} by using suitable changes of scales and Harnack-type inequalities. In Section~\ref{SecApriori2}, we prove that
\begin{equation}\label{Firstorder}
\left|\left|x\right|^\mu u\(x\)-\alpha\right|+\left|\left|x\right|^{\mu+1}\nabla u\(x\)+\alpha\mu\left|x\right|^{-1}x\right|\to0\quad\text{as }\left|x\right|\to\infty
\end{equation}
by using an approach based on comparison arguments. In Section~\ref{SecKelvin}, we then reduce the problem at infinity to another elliptic-type problem near the origin by using a Kelvin-type transformation. More precisely, we prove the following:

\begin{proposition}\label{Pr2}
Let $n$, $p$, $a$, $b$, $q$, $\mu$ and $f$ be as in Theorem~\ref{Th1}, $R_0>0$ and $u$ be a weak solution of the equation
\begin{equation}\label{Pr2Eq1}
-\divergence\(\left|x\right|^{-ap}\left|\nabla u\right|^{p-2}\nabla u\) =f\(x,u\)\quad\text{in }\R^n\backslash B\(0,R_0\).
\end{equation}
Assume that there exists a constant $\alpha>0$ such that \eqref{Firstorder} holds true. Let $r_0:=1/R_0$ and $u_\ast:B\(0,r_0\)\backslash\left\{0\right\}\to\R$ be the function defined by 
\begin{equation}\label{Pr2Eq2}
u_\ast\(y\):=\left|y\right|^{-\mu}u\(\left|y\right|^{-2}y\)\quad\forall y\in B\(0,r_0\).
\end{equation}
Then there exists $r\in\(0,r_0\)$ such that $u_\ast$ is a weak solution (see Definition~\ref{Def1} below) to an equation of the form
\begin{equation}\label{Pr2Eq3}
-\divergence\(A\(y,u_\ast,\nabla u_\ast\)\)=B\(y,u_\ast,\nabla u_\ast\)\quad\text{in }B\(0,r\)
\end{equation}
for some Caratheodory functions $A:B\(0,r\)\times\R\times\R^n\to\R^n$ and $B:B\(0,r\)\times\R\times\R^n\to\R$ satisfying
\begin{equation}\label{Pr2Eq4}
\left\{\begin{aligned}&\left|A\(y,z,\xi\)\right|\le C\left|y\right|^\gamma\left|\xi\right|\\&A\(y,z,\xi\)\cdot\xi\ge C^{-1}\left|y\right|^\gamma\left|\xi\right|^2\\&\left|B\(y,z,\xi\)\right|\le C\big(\left|y\right|^\gamma\left|\xi\right|^2+\left|y\right|^{\gamma'}\big),\end{aligned}\right.
\end{equation}
where 
\begin{equation}\label{Pr2Eq5}
\gamma:=\mu+2-n,\quad\gamma':=\frac{\mu q}{p}-n
\end{equation} 
and $C$ is a positive constant independent of $\(y,z,\xi\)\in B\(0,r\)\times\R\times\R^n$.
\end{proposition}

The definition of weak solution of \eqref{Pr2Eq3} is as follows:

\begin{definition}\label{Def1}
For every $r>0$ and $\gamma,\gamma'>-n$, we define the space $H^{1,2}\(B\(0,r\),\left|y\right|^{\gamma}\)$ as the completion of $C^\infty\(B\(0,r\)\)$ with respect to the norm 
$$\left\|u\right\|_{H^{1,2}\(B\(0,r\),\left|y\right|^{\gamma}\)}=\(\int_{B\(0,r\)}\left|y\right|^{\gamma}\left|u\right|^2dy+\int_{B\(0,r\)}\left|y\right|^{\gamma}\left|\nabla u\right|^2dy\)^{\frac12}.$$
Given two Caratheodory functions $A:B\(0,r\)\times\R\times\R^n\to\R^n$ and $B:B\(0,r\)\times\R\times\R^n\to\R$ satisfying \eqref{Pr2Eq4}, we say that $u$ is a weak solution of the equation
$$-\divergence\(A\(y,u,\nabla u\)\)=B\(y,u,\nabla u\)\quad\text{in }B\(0,r\)$$
if $u\in H^{1,2}\(B\(0,r\),\left|y\right|^{\gamma}\)$ and 
$$\int_{B\(0,r\)}\(A\(y,u,\nabla u\)\cdot\nabla\varphi-B\(y,u,\nabla u\)\varphi\)dy=0\,\,\,\,\,\forall\varphi\in C_0^\infty\(B\(0,r\)\).$$ 
\end{definition}

It is well known that in the case where $f\(x,u\)=\Lambda\left|x\right|^{-bq}\left|u\right|^{q-2}u$ and $p=2$, the problem \eqref{Main Problem} is invariant under the transformation $u\mapsto u_\ast$. With Proposition~\ref{Pr2}, we show that a similar transformation can still be applied in the general case provided we first establish \eqref{Firstorder}. It is interesting to remark than even in the unweighted case where $a=b=0$, the problem that we obtain after transformation is actually of weighted-type when $p\ne2$. At the end of Section~\ref{SecKelvin}, we then combine Proposition~\ref{Pr2} with an H\"older continuity result of Stredulinski~\cite{Str} to obtain that \eqref{Th1Eq4} holds true, thus completing the proof of Theorem~\ref{Th1}. For other references on the H\"older continuity of solutions to weighted elliptic equations, let us mention for instance the works of Colorado and Peral~\cite{ColPer}, Di Fazio and  Zamboni~\cite{DiFZam}, Felli and Schneider~\cite{FelSch2} and Monticelli, Rodney and Wheeden~\cite{MonRodWhe2}.  

\begin{remark}
It is possible to prove a more general version of Proposition~\ref{Pr2} in the case where $u_\ast$ is defined as 
$$u_\ast\(y\):=\left|y\right|^{-\mu}u\(\left|y\right|^{-\sigma}y\)\quad\forall y\in B\(0,r_0\)$$
for some constant $\sigma>1$. In this case, we obtain that $u_\ast$ solves an equation of the form \eqref{Pr2Eq3} for some Caratheodory functions $A:B\(0,r\)\times\R\times\R^n\to\R^n$ and $B:B\(0,r\)\times\R\times\R^n\to\R$ satisfying \eqref{Pr2Eq4} with 
$$\gamma:=\(\sigma-1\)\mu+2-n\quad\text{and}\quad\gamma':=\(\sigma-1\)\mu\frac{q}{p}-n.$$
In particular, when $n\ge3$, we can choose 
$$\sigma:=1+\frac{n-2}{\mu}\,,$$ 
which gives $\gamma=0$. Furthermore, in this case, a straightforward computation gives 
$$\gamma'=\(n-2\)\frac{q}{p}-n,$$
 which is greater than $-2$ since $q>p$. In particular, this allows to apply the H\"older continuity result of Trudinger~\cite{Tru1}*{Theorem~5.2} to obtain that \eqref{Th1Eq4} holds true.
\end{remark}

\section{Regularity and boundedness results}\label{SecRegularity}

We start by introducing some suitable function spaces.

\begin{definition}
For every $s>0$, $\gamma\in\R$ and measurable set $\Omega \subseteq \R^n$, we define $L^s\(\Omega,\left|x\right|^\gamma\)$ as the set of all measurable functions $u:\Omega\to\R$ such that  
$$\|u\|_{L^s\(\Omega,\left|x\right|^\gamma\)}:= \( \int_{\Omega}\left|x\right|^\gamma\left|u\right|^s dx   \)^{\frac1s} < \infty.$$
Furthermore, we denote $L^s\(\Omega\):=L^s\(\Omega, 1\)$.
\end{definition}

\begin{definition}
For every $s>0$, $\gamma\in\R$ and measurable set $\Omega\subseteq\R^n$, we define $ L^{s,\infty}\( \Omega,\left|x\right|^\gamma\)$ as the set of all measurable functions $u:\Omega\to\R$ such that 
$$\left\|u\right\|_{L^{s,\infty}\(\Omega,\left|x\right|^\gamma\)}:=\sup_{h>0} \( h \(\int_{W_h}\left|x\right|^\gamma dx\)^{\frac{1}{s}}\)<\infty,$$
where 
\begin{equation}\label{Level}
W_h:=\left\{x\in\R^n:\,\left|u\(x\)\right|>h\right\}.
\end{equation}
Furthermore, we denote $L^{s,\infty}\(\Omega\):=L^{s,\infty}\(\Omega, 1\)$.
\end{definition}

As a first step in the proof of Theorem~\ref{Th1}, we prove the following:

\begin{step}\label{Step:boundedness for solutions}
Let $n$, $p$, $a$, $b$, $q$ and $f$ be as in Theorem~\ref{Th1}. Then every solution of \eqref{Main Problem} belongs to $L^\infty\(\R^n\)$.
\end{step}

\proof
Here we adapt some ideas which originates from Trudinger~\cite{Tru2}. Let $u$ be a solution of \eqref{Main Problem}. We begin with proving that $\left|u\right|^{\beta/p}\in L^q\big(\R^n,\left|x\right|^{-bq}\big)$ for all $\beta>p$. For every $h>0$, we define
$$\varphi_h\(u\):=\min\big(\left|u\right|,h\big)^{\frac{\beta-p}{p}}.$$
By using $\varphi_h\(u\)^pu$ as a test function, we obtain
\begin{multline}\label{PrEq1}
\int_{\R^n}\left|x\right|^{-ap}\left|\nabla u\right|^p\big(h^{\beta-p}\chi_{W_h}+\(\beta-p+1\)\left|u\right|^{\beta-p}\chi_{\R^n\backslash W_h}\big)dx\\
=\int_{\R^n}f\(x,u\)\varphi_h\(u\)^pu\,dx,
\end{multline}
where $W_h$ is as in \eqref{Level} and $\chi_{W_h}$ and $\chi_{\R^n\backslash W_h}$ are the characteristic functions of the sets $W_h$ and $\R^n\backslash W_h$, respectively. Since $\beta>p$, it follows from \eqref{Th1Eq1} and \eqref{PrEq1} that
\begin{equation}\label{PrEq2}
\int_{\R^n}\left|x\right|^{-ap}\left|\nabla u\right|^p\varphi_h\(u\)^pdx\le\Lambda\int_{\R^n}\left|x\right|^{-bq}\left|u\right|^q\varphi_h\(u\)^pdx.
\end{equation}
On the other hand, for every $k\in\(0,h\)$, we have
\begin{align}\label{PrEq3}
&\int_{\R^n}\left|x\right|^{-bq}\left|u\right|^q\varphi_h\(u\)^pdx\nonumber\\
&\qquad\le k^{\beta-p}\int_{\R^n}\left|x\right|^{-bq}\left|u\right|^qdx+\int_{W_k}\left|x\right|^{-bq}\left|u\right|^q\varphi_h\(u\)^pdx.
\end{align}
Since $q>p$, by using H\"older's inequality, we obtain
\begin{align}\label{PrEq4}
&\int_{W_k}\left|x\right|^{-bq}\left|u\right|^q\varphi_h\(u\)^pdx\nonumber\\
&\qquad\le\(\int_{W_k}\left|x\right|^{-bq}\left|u\right|^qdx\)^{\frac{q-p}{q}}\(\int_{W_k}\left|x\right|^{-bq}\left|u\right|^q\varphi_h\(u\)^qdx\)^{\frac{p}{q}}.
\end{align}
By applying the Caffarelli--Kohn--Nirenberg inequality to the function $\varphi_h\(u\)u$, we obtain
\begin{align}\label{PrEq5}
&\(\int_{\R^n}\left|x\right|^{-bq}\left|u\right|^q\varphi_h\(u\)^qdx\)^{\frac{p}{q}}\nonumber\\
&\qquad\le K\int_{\R^n}\left|x\right|^{-ap}\left|\nabla u\right|^p\big(h^{\beta-p}\chi_{W_h}+\(\beta/p\)^p\left|u\right|^{\beta-p}\chi_{\R^n\backslash W_h}\big)dx\nonumber\\
&\qquad\le K\(\beta/p\)^p\int_{\R^n}\left|x\right|^{-ap}\left|\nabla u\right|^p\varphi_h\(u\)^pdx
\end{align}
for some constant $K=K\(n,a,b,p\)>0$. It follows from \eqref{PrEq2}--\eqref{PrEq5} that
\begin{multline}\label{PrEq7}
\int_{\R^n}\left|x\right|^{-ap}\left|\nabla u\right|^p\varphi_h\(u\)^pdx\le\Lambda\Bigg(k^{\beta-p}\int_{\R^n}\left|x\right|^{-bq}\left|u\right|^qdx\\
+K\(\beta/p\)^p\(\int_{W_k}\left|x\right|^{-bq}\left|u\right|^qdx\)^{\frac{q-p}{q}}\hspace{-2pt}\int_{\R^n}\left|x\right|^{-ap}\left|\nabla u\right|^p\varphi_h\(u\)^pdx\Bigg).
\end{multline}
Since $u\in L^q\big(\R^n,\left|x\right|^{-bq}\big)$ and $q>p$, we have
$$\lim_{k\to+\infty}\(\int_{W_k}\left|x\right|^{-bq}\left|u\right|^qdx\)^{\frac{q-p}{q}}=0.$$
Therefore, by choosing $k$ sufficiently large (depending on $u$), it follows from \eqref{PrEq7} that
$$\int_{\R^n}\left|x\right|^{-ap}\left|\nabla u\right|^p\varphi_h\(u\)^pdx\le2\Lambda k^{\beta-p}\int_{\R^n}\left|x\right|^{-bq}\left|u\right|^qdx,$$
which together with \eqref{PrEq5} gives
\begin{equation}\label{PrEq8}
\(\int_{\R^n}\left|x\right|^{-bq}\left|u\right|^q\varphi_h\(u\)^qdx\)^{\frac{p}{q}}\le 2\Lambda k^{\beta-p}K\(\beta/p\)^p\int_{\R^n}\left|x\right|^{-bq}\left|u\right|^qdx.
\end{equation}
Since $u\in L^q\big(\R^n,\left|x\right|^{-bq}\big)$, by passing to the limit as $h\to+\infty$ into \eqref{PrEq8}, we then obtain that $\left|u\right|^{\beta/p}\in L^q\big(\R^n,\left|x\right|^{-bq}\big)$. Now we prove that $u\in L^\infty\(\R^n\)$. For every $h>0$, we define
$$\psi_h\(u\):=\sgn\(u\)\max\big(\left|u\right|-h,0\big),$$
where $\sgn\(u\)$ denotes the sign of $u$. By using $\psi_h\(u\)$ as a test function, we obtain
\begin{equation}\label{PrEq9}
\int_{W_h}\left|x\right|^{-ap}\left|\nabla u\right|^pdx=\int_{W_h}f\(x,u\)\psi_h\(u\)dx,
\end{equation}
where $W_h$ is as in \eqref{Level}. For every $\beta>p$, it follows from \eqref{Th1Eq1}, \eqref{PrEq9} and H\"older's inequality that
\begin{multline}\label{PrEq10}
\int_{W_h}\left|x\right|^{-ap}\left|\nabla u\right|^pdx\le\Lambda\(\int_{W_h}\left|x\right|^{-bq}dx\)^{1-\frac{1}{q}-\frac{p\(q-1\)}{\beta q}}\\
\times\(\int_{W_h}\left|x\right|^{-bq}\left|u\right|^{\frac{\beta q}{p}}dx\)^{\frac{p\(q-1\)}{\beta q}}\(\int_{W_h}\left|x\right|^{-bq}\left|\psi_h\(u\)\right|^qdx\)^{\frac{1}{q}}.
\end{multline}
By applying the Caffarelli--Kohn--Nirenberg inequality, we obtain
\begin{equation}\label{PrEq11}
\(\int_{W_h}\left|x\right|^{-bq}\left|\psi_h\(u\)\right|^qdx\)^{\frac{p}{q}}\le K\int_{W_h}\left|x\right|^{-ap}\left|\nabla u\right|^pdx.
\end{equation}
Since $\left|u\right|^{\beta/p}\in L^q\big(\R^n,\left|x\right|^{-bq}\big)$, it follows from \eqref{PrEq10} and \eqref{PrEq11} that
\begin{equation}\label{PrEq12}
\int_{W_h}\left|x\right|^{-bq}\left|\psi_h\(u\)\right|^qdx\le C\(\int_{W_h}\left|x\right|^{-bq}dx\)^{\frac{q}{p-1}\(1-\frac{1}{q}-\frac{p\(q-1\)}{\beta q}\)}
\end{equation}
for some constant $C>0$ independent of $h$. It follows from H\"older's inequality and \eqref{PrEq12} that
\begin{align}\label{PrEq13}
&\int_{W_h}\left|x\right|^{-bq}\left|\psi_h\(u\)\right|dx\nonumber\\
&\qquad\le\(\int_{W_h}\left|x\right|^{-bq}dx\)^{\frac{q-1}{q}}\(\int_{W_h}\left|x\right|^{-bq}\left|\psi_h\(u\)\right|^qdx\)^{\frac{1}{q}}\nonumber\\
&\qquad\le C^{\frac{1}{q}}\(\int_{W_h}\left|x\right|^{-bq}dx\)^{\frac{q-1}{q}+\frac{1}{p-1}\(1-\frac{1}{q}-\frac{p\(q-1\)}{\beta q}\)}.
\end{align}
On the other hand, by applying Tonelli's theorem, we obtain
\begin{equation}\label{PrEq14}
\int_{W_h}\left|x\right|^{-bq}\left|\psi_h\(u\)\right|dx=\int_h^\infty\int_{W_s}\left|x\right|^{-bq}dx\,ds.
\end{equation}
Furthermore, by choosing $\beta$ large enough so that $\beta>p\(q-1\)/\(q-p\)$, we obtain
\begin{equation}\label{PrEq15}
\frac{q-1}{q}+\frac{1}{p-1}\(1-\frac{1}{q}-\frac{p\(q-1\)}{\beta q}\)>1.
\end{equation}
Therefore, since $\left|x\right|^{-bq}\in L^1\(\R^n\)$, it follows from \eqref{PrEq13}--\eqref{PrEq15} that
$$\int_{W_h}\left|x\right|^{-bq}dx=0\quad\text{for large }h>0$$
and so $u\in L^\infty\(\R^n\)$. This ends the proof of Step~\ref{Step:boundedness for solutions}.
\endproof

We then obtain the following:

\begin{step}\label{Step:Holder continuity}
Let $n$, $p$, $a$, $b$, $q$ and $f$ be as in Theorem~\ref{Th1}. Then every solution of \eqref{Main Problem} belongs to $C^1\(\R^n\backslash\left\{0\right\}\)$.
\end{step}

\proof[Proof]
This step follows directly from Step~\ref{Step:boundedness for solutions} together with the regularity results of DiBenedetto~\cite{DiB} and Tolksdorf~\cite{Tol} (see also Evans~\cite{Eva}, Lewis~\cite{Lew}, Uhlenbeck~\cite{Uhl} and Ural{\cprime}ceva~\cite{Ura} for previous results on this question).
\endproof

%

The next result is concerned with the boundedness of solutions of \eqref{Main Problem} in weak Lebesgue spaces. 

\begin{step}\label{Step:solutions are L^(p_*-1)}
Let $n$, $p$, $a$, $b$, $q$ and $f$ be as in Theorem~\ref{Th1}. Then every solution of \eqref{Main Problem} belongs to $L^{q-q/p, \infty}\big(\R^n, \left|x\right|^{-bq}\big)$.
\end{step}

\begin{proof}
Let $u$ be a non-trivial solution of \eqref{Main Problem}. For every $h>0,$ we define 
$$\tau_h\(u\):= \sgn\(u\) \min\(\left|u\right|,h\),$$ 
where $\sgn\(u\)$ denotes the sign of $u$. By testing (\ref{Main Problem}) with $\tau_h\(u\)$ and using \eqref{Th1Eq1}, we obtain 
\begin{multline}\label{formula:testing by tau_h(u)}
\int_{\R^n\backslash W_h} \left|x\right|^{-ap} \left|\nabla u\right|^p dx \le\Lambda\bigg(\int_{\R^n\backslash W_h} \left|x\right|^{-bq} \left|u\right|^q dx\\
+ h \int_{W_h} \left|x\right|^{-bq} \left|u\right|^{q-1} dx\bigg),
\end{multline}
where $W_h$ is as in \eqref{Level}. On the other hand, straightforward computations give
\begin{equation}\label{Equality: expression of C-K-N term using tau_h}
\int_{\R^n\backslash W_h} \left|x\right|^{-bq} \left|u\right|^q dx = \int_{\R^n}\left|x\right|^{-bq}\left|\tau_h\(u\)\right|^{q} dx- h^q \int_{W_h} \left|x\right|^{-bq}dx 
\end{equation}
and 
\begin{multline}\label{Equality: expression of C-K-N term on |u|>h}
\int_{W_h} \left|x\right|^{-bq} \left|u\right|^{q-1} dx =  \(q-1\) \int_h^{\infty} \left( \int_{W_s}  \left|x\right|^{-bq} dx \right) s^{q-2} ds\\
+ h^{q-1}  \int_{W_h}  \left|x\right|^{-bq} dx.
\end{multline}
Plugging \eqref{Equality: expression of C-K-N term using tau_h} and \eqref{Equality: expression of C-K-N term on |u|>h} into (\ref{formula:testing by tau_h(u)}), we obtain 
\begin{multline}\label{Equality: expression of energy C-K-N term on |u|<h}
\int_{\R^n\backslash W_h} \left|x\right|^{-ap} \left|\nabla u\right|^p dx\le\Lambda\bigg(\int_{\R^n}\left|x\right|^{-bq}\left|\tau_h\(u\)\right|^{q} dx\\
+\(q-1\)h\int_h^{\infty} \left( \int_{W_s}  \left|x\right|^{-bq} dx \right) s^{q-2} ds\bigg).
\end{multline}
The Caffarelli--Kohn--Nirenberg inequality yields
\begin{equation}\label{Inequality:C-K-N for tau_h on |u| <h}
  \int_{\R^n}\left|x\right|^{-bq}\left|\tau_h\(u\)\right|^{q} dx  \le K \left(  \int_{\R^n\backslash W_h} \left|x\right|^{-ap} \left|\nabla u\right|^p dx \right)^{\frac{q}{p}}
\end{equation}
for some constant $K=K\(n,a,b,p\)>0$. Since $u\in D^{1,p}\(\R^n,\left|x\right|^{-ap}\)$, it follows from (\ref{Equality: expression of C-K-N term using tau_h}), (\ref{Equality: expression of energy C-K-N term on |u|<h}) and (\ref{Inequality:C-K-N for tau_h on |u| <h}) that there exists a constant $C=>0$ such that
\begin{multline}\label{Inequality:estimate fo G(h)}
h^q \int_{W_h} \left|x\right|^{-bq}dx \le  \int_{\R^n}\left|x\right|^{-bq}\left|\tau_h\(u\)\right|^{q} dx\\
\le C \left( h  \int_h^{\infty} \left( \int_{W_s}  \left|x\right|^{-bq} dx \right) s^{q-2} ds \right)^{\frac{q}{p}} \quad \text{for small } h>0.
\end{multline}
We now define 
$$G\(h\) := \left( \int_h^\infty g\(s\) ds \right)^{\frac{p-q}{p}},\quad\text{where}\quad g\(s\) := s^{q-2} \int_{W_s}  \left|x\right|^{-bq} dx. $$
In particular, $G$ is positive, non-decreasing and locally absolutely continuous in $\big(0,\left\|u\right\|_{L^\infty\(\R^n\)}\big)$ with derivative
$$G'(h) =\frac{q-p}{p}\left( \int_h^\infty g\(s\) ds\right)^{-\frac{q}{p}} g\(h\) \quad \text{ for a.e. } h \in (0, \|u\|_{L^\infty(\R^n)}).$$
By using \eqref{Inequality:estimate fo G(h)}, we then obtain
\begin{equation}\label{G'(h)}
G'\(h\) \le\frac{q-p}{p}C h^{\frac{q-2p}{p}} \quad \text{for small } h>0.
\end{equation}
By integrating \eqref{G'(h)}, we obtain 
\begin{equation}\label{Inequality:upper bound for G(h)-G(0)}
G\(h\) - G\(0\) \le  C h^{\frac{q-p}{p}} \quad \text{for small }h>0,
\end{equation}
where $G\(0\)$ stands for the limit of $G\(h\)$ as $h\to0$.
On the other hand, by using (\ref{Equality: expression of C-K-N term on |u|>h}) together with dominated convergence, we obtain
$$\(q-1\)hG\(h\)^{\frac{-p}{q-p}}\le h\int_{W_h}\left|x\right|^{-bq}\left|u\right|^{q-1}dx=\smallo\(1\) \quad \text{as } h\to0.$$
This coupled with (\ref{Inequality:upper bound for G(h)-G(0)}) yields that $G\(0\)>0$. By using \eqref{Inequality:estimate fo G(h)} and since  $G$ is non-decreasing, we then obtain 
$$h^{\frac{q\(p-1\)}{p}}
\int_{W_h}  \left|x\right|^{-bq} dx  \le C G\(h\)^{-\frac{q}{q-p}}\le C G\(0\)^{-\frac{q}{q-p}} \quad \text{for small } h>0, $$
which implies that $u\in L^{q-q/p,\infty}\big(\R^n,\left|x\right|^{-bq}\big)$. This ends the proof of Step~\ref{Step:solutions are L^(p_*-1)}.
\end{proof}

\section{The upper bound estimates}\label{SecApriori1}

This section is devoted to the proof of \eqref{Th1Eq3}. We begin with establishing a decay estimate, which is weaker than \eqref{Th1Eq3}, but which will serve as a preliminary step in the proof of \eqref{Th1Eq3}. 

\begin{step}\label{Th1Step1}
Let $n$, $p$, $a$, $b$, $q$ and $f$ be as in Theorem~\ref{Th1}. Let $u$ be a solution of \eqref{Main Problem}. Then there exists a constant $K_0>0$ such that
\begin{equation}\label{Th1Step1Eq1}
\left|x\right|^{\nu}\left|u\(x\)\right|+\left|x\right|^{\nu+1}\left|\nabla u\(x\)\right|\le K_0\quad\forall x\in\R^n\backslash B\(0,1/2\),
\end{equation}
where $\nu:=\(n-p\(1+a\)\)/p$.
\end{step}

\proof
For every $R>0$, we define
$$u_R\(x\):=R^{\nu}u\(Rx\)\quad\forall x\in\R^n.$$
As is easily seen, in order to prove \eqref{Th1Step1Eq1}, it suffices to show that
\begin{equation}\label{Th1Step1Eq2}
\left|u_R\(x\)\right|+\left|\nabla u_R\(x\)\right|\le K_0\quad\forall R>1/2,\,x\in\S^n
\end{equation} 
for some constant $K_0>0$. By using \eqref{Main Problem}, we obtain 
\begin{align}\label{Th1Step1Eq3}
&-\Delta_p u_R=R^{\nu\(p-1\)+p}\big(\left|Rx\right|^{ap}f\(Rx,u\(Rx\)\)\nonumber\\
&\quad-ap\left|Rx\right|^{-2}\left|\nabla u\(Rx\)\right|^{p-2}\<\nabla u\(Rx\),Rx\>\big)\nonumber\\
&=R^{\nu\(p-1\)+\(a+1\)p}\left|x\right|^{ap}f\big(Rx,R^{-\nu}u_R\big)\nonumber\\
&-ap\left|x\right|^{-2}\left|\nabla u_R\right|^{p-2}\<\nabla u_R,x\>=:g_R\(x,u_R,\nabla u_R\)\quad\text{in }\R^n\backslash\left\{0\right\}.
\end{align}
By using \eqref{Th1Eq1} together with straightforward computations, we obtain
\begin{multline}\label{Th1Step1Eq4}
\left|g_R\(x,u_R,\nabla u_R\)\right|\le\Lambda R^{\nu\(p-q\)+\(a+1\)p-bq}\left|x\right|^{ap-bq}\left|u_R\right|^{q-1}\\
+\left|a\right|p\left|x\right|^{-1}\left|\nabla u_R\right|^{p-1}=\Lambda \left|x\right|^{ap-bq}\left|u_R\right|^{q-1}+\left|a\right|p\left|x\right|^{-1}\left|\nabla u_R\right|^{p-1}
\end{multline}
for a.e. $x\in\R^n$. Assuming by contradiction that \eqref{Th1Step1Eq2} does not hold true, we obtain that there exist sequences $\(R_k\)_{k\in\N}$ in $\(1/2,\infty\)$ and $\(y_k\)_{k\in\N}$ in $\S^n$ such that 
\begin{equation}\label{Th1Step1Eq2Temp}
\left|u_{R_k}\(y_k\)\right|+\left|\nabla u_{R_k}\(y_k\)\right|\ge k\quad\forall k\in\R^n.
\end{equation} 
By using a doubling property (see Pol\'a\v{c}ik, Quittner and Souplet~\cite{PolQuiSou}*{Lemma~5.1}), we then obtain that there exists a sequence $\(x_k\)_{k\in\N}$ in $B\(0,2\)\backslash B\(0,1/2\)$ such that 
\begin{equation}\label{Th1Step1Eq5}
\lambda_k:=\Big(\left|u_{R_k}\(x_k\)\right|^{\frac{q}{n}}+\left|\nabla u_{R_k}\(x_k\)\right|^{\frac{q}{n+q}}\Big)^{-1}\to0\quad\text{as }k\to\infty,
\end{equation}
\begin{equation}\label{Th1Step1Eq6}
B\(x_k,2k\lambda_k\)\subset B\(0,2\)\backslash B\(0,1/2\)\quad\forall k\in\N
\end{equation}
and
\begin{equation}\label{Th1Step1Eq7}
\left|\widetilde{u}_k\(y\)\right|^{\frac{q}{n}}+\left|\nabla \widetilde{u}_k\(y\)\right|^{\frac{q}{n+q}}\le2\quad\forall k\in\N,\, y\in B\(0,k\),
\end{equation}
where
$$\widetilde{u}_k\(y\):=\lambda_k^{\frac{n}{q}}u_{R_k}\(x_k+\lambda_ky\).$$
By using \eqref{Th1Step1Eq3}, we obtain
\begin{align*}
-\Delta_p \widetilde{u}_k&=\lambda_k^{\frac{n}{q}\(p-1\)+p}g_{R_k}\(x_k+\lambda_ky,u_{R_k}\(x_k+\lambda_ky\),\nabla u_{R_k}\(x_k+\lambda_ky\)\)\\
&=\lambda_k^{\frac{n}{q}\(p-1\)+p}g_{R_k}\Big(x_k+\lambda_ky,\lambda_k^{-\frac{n}{q}}\widetilde{u}_k,\lambda_k^{-\frac{n+q}{q}}\nabla \widetilde{u}_k\Big)
\end{align*}
in $\R^n\backslash\left\{-x_k/\lambda_k\right\}$. Furthermore, by using \eqref{Th1Step1Eq4}--\eqref{Th1Step1Eq7} and observing that $\frac{n}{q}\(p-q\)+p=p\(b-a\)\ge0$, we obtain
\begin{multline*}
\Big|\lambda_k^{\frac{n}{q}\(p-1\)+p}g_{R_k}\Big(x_k+\lambda_ky,\lambda_k^{-\frac{n}{q}}\widetilde{u}_k,\lambda_k^{-\frac{n+q}{q}}\nabla \widetilde{u}_k\Big)\Big|\\
\le\Lambda\lambda_k^{p\(b-a\)}\left|x_k+\lambda_ky\right|^{ap-bq}\left|\widetilde{u}_k\right|^{q-1}+\lambda_k\left|a\right|p\left|x_k+\lambda_ky\right|^{-1}\left|\nabla \widetilde{u}_k\right|^{p-1}\le C
\end{multline*}
for a.e. $y\in B\(0,k\)$, for some constant $C>0$ independent of $k$. Since $\left|x_k\right|/\lambda_k\to\infty$ as $k\to\infty$, by applying the results of DiBenedetto~\cite{DiB} and Tolksdorf~\cite{Tol}, we then obtain that $\(\widetilde{u}_k\)_{k\in\N}$ is bounded in $C^{1,\theta}_{\loc}\(\R^n\)$ for some $\theta\in\(0,1\)$ and so $\(\widetilde{u}_k\)_{k\in\N}$ converges up to a subsequence in $C^1_{\loc}\(\R^n\)$ to some function $\widetilde{u}_\infty\in C^1\(\R^n\)$. Remark that by definition of $\widetilde{u}_k$, we have
\begin{equation}\label{Th1Step1Eq8}
\left|\widetilde{u}_k\(0\)\right|^{\frac{q}{n}}+\left|\nabla \widetilde{u}_k\(0\)\right|^{\frac{q}{n+q}}=1.
\end{equation}
By passing to the limit as $k\to\infty$ into \eqref{Th1Step1Eq8}, we then obtain
\begin{equation}\label{Th1Step1Eq9}
\left|\widetilde{u}_\infty\(0\)\right|^{\frac{q}{n}}+\left|\nabla \widetilde{u}_\infty\(0\)\right|^{\frac{q}{n+q}}=1.
\end{equation}
On the other hand, by using \eqref{Th1Step1Eq6} together with straightforward computations, we obtain
\begin{align}\label{Th1Step1Eq10}
\int_{B\(0,k\)}\left|\widetilde{u}_k\right|^qdy&=\int_{B\(x_k,k\lambda_k\)}\left|u_{R_k}\right|^qdx\le\int_{B\(0,2\)\backslash B\(0,1/2\)}\left|u_{R_k}\right|^qdx\nonumber\\
&\le2^{\left|b\right|q}\int_{B\(0,2\)\backslash B\(0,1/2\)}\left|x\right|^{-bq}\left|u_{R_k}\right|^qdx\nonumber\\
&\le2^{\left|b\right|q}\int_{B\(0,2R_k\)\backslash B\(0,R_k/2\)}\left|x\right|^{-bq}\left|u\right|^qdx.
\end{align}
Since $R_k>1/2$ and $u\in C^1\(\R^n\backslash\left\{0\right\}\)$, it follows from \eqref{Th1Step1Eq5} and \eqref{Th1Step1Eq6} that $R_k\to\infty$ as $k\to\infty$. Since $u\in L^q\big(\R^n,\left|x\right|^{-bq}\big)$, it then follows from \eqref{Th1Step1Eq10} that $\widetilde{u}_k\to0$ in $L^q_{\loc}\(\R^n\)$ and so $\widetilde{u}_\infty\equiv0$ in $\R^n$. This is in contradiction with \eqref{Th1Step1Eq9}. This ends the proof of Step~\ref{Th1Step1}.
\endproof

The next step is as follows:

\begin{step}\label{Th1Step2}
Let $n$, $p$, $a$, $b$, $q$ and $f$ be as in Theorem~\ref{Th1}. Let $u$ be a solution of \eqref{Main Problem}. Let $K_0$ be as in Step~\ref{Th1Step1}. Then $v:=\left|u\right|$ satisfies the inequality
\begin{equation}\label{Th1Step2Eq1}
-\Delta_pv\le \overline{g}\(x,v,\nabla v\)\quad\text{in }\R^n\backslash B\(0,1/2\),
\end{equation}
where
\begin{equation}\label{Th1Step2Eq2}
\overline{g}\(x,v,\nabla v\):=\Lambda K_0^{q-p}\left|x\right|^{-p}v^{p-1}+\left|a\right|p\left|x\right|^{-1}\left|\nabla v\right|^{p-1}.
\end{equation}
\end{step}

The inequality \eqref{Th1Step2Eq1} is to be understood in the sense that
$$\int_{\R^n}\left|\nabla v\right|^{p-2}\<\nabla v,\nabla\varphi\>dx\le\int_{\R^n}\overline{g}\(x,v,\nabla v\)\varphi\,dx$$
for all $\varphi\in C^\infty_c\(\R^n\backslash B\(0,1/2\)\)$ such that $\varphi\ge0$ in $\R^n\backslash B\(0,1/2\)$.

\proof
As is easily seen, the equation satisfied by $u$ can be rewritten as
\begin{equation}\label{Th1Step2Eq3}
-\Delta_pu=\left|x\right|^{ap}f\(x,u\)-ap\left|x\right|^{-2}\left|\nabla u\right|^{p-2}\<\nabla u,x\>\quad\text{in }\R^n\backslash\left\{0\right\}.
\end{equation}
By using \eqref{Th1Eq1} and \eqref{Th1Step1Eq1} together with straightforward computations, we obtain
\begin{equation}\label{Th1Step2Eq4}
\left|\left|x\right|^{ap}f\(x,u\)-ap\left|x\right|^{-2}\left|\nabla u\right|^{p-2}\<\nabla u,x\>\right|\le \overline{g}\(x,v,\nabla v\)
\end{equation}
for a.e. $x\in B\(0,1/2\)$, where $\overline{g}\(x,v,\nabla v\)$ is as in \eqref{Th1Step2Eq2}. The inequality \eqref{Th1Step2Eq1} then follows from \eqref{Th1Step2Eq3} and  \eqref{Th1Step2Eq4} by applying an extended version of Kato's inequality~\cite{Kato} for the $p$--Laplace operator (see Cuesta~Leon~\cite{Cue}*{Proposition~3.2}). This ends the proof of Step~\ref{Th1Step2}.
\endproof

We can now prove \eqref{Th1Eq3} by using Steps~\ref{Th1Step1} and~\ref{Th1Step2}.

\proof[Proof of \eqref{Th1Eq3}]
Let $u$ be a solution of \eqref{Main Problem} and $v:=\left|u\right|$. For every $R>0$, we define
\begin{align}\label{Th1Step3Eq2}
u_R\(x\):=R^\mu u\(Rx\)\quad\text{and}\quad v_R\(x\):=R^\mu v\(Rx\)\quad\forall x\in\R^n,
\end{align}  
where $\mu$ is as in \eqref{Th1Eq4}. As is easily seen, in order to prove \eqref{Th1Eq3}, it suffices to show that
\begin{equation}\label{Th1Step3Eq3}
v_R\(x\)+\left|\nabla u_R\(x\)\right|\le C_0\quad\forall R>1,\,x\in\S^n
\end{equation} 
for some constant $C_0>0$. By using \eqref{Th1Step2Eq1} and remarking that
$$R^{n-pa}\overline{g}\big(Rx,v\(Rx\),\nabla v\(Rx\)\big)=\overline{g}\(x,v_R\(x\),\nabla v_R\(x\)\),$$
we obtain that $v_R$ satisfies the inequality
$$-\Delta_pv_R\le \overline{g}\(x,v_R,\nabla v_R\)\quad\text{in }\R^n\backslash B\(0,1/\(2R\)\).$$
We can then apply a weak Harnack inequality (see Trudinger~\cite{Tru1}*{Theorem~1.3}), which gives that for every $\epsilon>0$, there exists a constant $c_\epsilon>0$ such that
\begin{equation}\label{Th1Step3Eq4}
\left\|v_R\right\|_{L^\infty\(B\(0,2\)\backslash B\(0,1/2\)\)}\le c_\epsilon\left\|v_R\right\|_{L^{p-1+\epsilon}\(B\(0,3\)\backslash B\(0,1/3\)\)}
\end{equation}
for all $R>1$. Since $q-q/p>p-1$, we can choose $\epsilon$ so that $p-1+\epsilon<q-q/p$. We then obtain
\begin{equation}\label{Th1Step3Eq5}
\left\|v_R\right\|_{L^{p-1+\epsilon}\(B\(0,3\)\backslash B\(0,1/3\)\)}\le c\left\|v_R\right\|_{L^{q-q/p,\infty}\(B\(0,3\)\backslash B\(0,1/3\)\)}
\end{equation}
for some constant $c=c\(n,p,a,b\)>0$. Furthermore, straightforward computations give
\begin{align}\label{Th1Step3Eq6}
\left\|v_R\right\|^{q-q/p}_{L^{q-q/p,\infty}\(B\(0,3\)\backslash B\(0,1/3\)\)}&=R^{-bq}\left\|u\right\|^{q-q/p}_{L^{q-q/p,\infty}\(B\(0,3R\)\backslash B\(0,R/3\)\)}\nonumber\\
&\le3^{\left|b\right|q}\left\|u\right\|^{q-q/p}_{L^{q-q/p,\infty}(\R^n,\left|x\right|^{-bq})}.
\end{align}
Since $u\in L^{q-q/p,\infty}\big(\R^n,\left|x\right|^{-bq}\big)$, it follows from \eqref{Th1Step3Eq4}--\eqref{Th1Step3Eq6} that
\begin{equation}\label{Th1Step3Eq7}
\left\|v_R\right\|_{L^\infty\(B\(0,2\)\backslash B\(0,1/2\)\)}\le C
\end{equation}
for some constant $C>0$ independent of $R>1$. We then infer \eqref{Th1Step3Eq3} from \eqref{Th1Step3Eq7} by applying the gradient estimates of DiBenedetto~\cite{DiB} and Tolksdorf~\cite{Tol}. This ends the proof of \eqref{Th1Eq3}.
\endproof

\section{The first-order term}\label{SecApriori2}

This section is devoted to the proof of \eqref{Firstorder}. We begin with proving the following:

\begin{step}\label{Th1Step4}
Let $n$, $p$, $a$, $b$, $q$, $\mu$ and $f$ be as in Theorem~\ref{Th1}. Let $u$ be a solution of \eqref{Main Problem} such that $u>0$ and $f\(x,u\)\ge0$ in $\R^n$. Then there exists a constant $c_0>0$ such that
\begin{equation}\label{Th1Step4Eq1}
u\(x\)\ge c_0\left|x\right|^{-\mu}\quad\forall x\in\R^n\backslash B\(0,1\).
\end{equation}
\end{step}

\proof
For every $R>0$, we let $u_R$ be the function defined as in \eqref{Th1Step3Eq2}. As is easily seen, in order to prove \eqref{Th1Step4Eq1}, it suffices to show that
\begin{equation}\label{Th1Step4Eq2}
u_R\(x\)\ge c_0\quad\forall R>1,\,x\in\S^n
\end{equation}
for some constant $c_0>0$. By using \eqref{Main Problem}, we obtain
\begin{equation}\label{Th1Step4Eq3}
-\divergence\(\left|x\right|^{-ap}\left|\nabla u_R\right|^{p-2}\nabla u_R\)=R^nf\big(Rx,R^{-\mu}u_R\big)\quad\text{in }\R^n.
\end{equation}
Furthermore, by using \eqref{Th1Eq1} and \eqref{Th1Eq3}, we obtain
\begin{equation}\label{Th1Step4Eq4}
\left|x\right|^{\mu}u_R\(x\)+\left|x\right|^{\mu+1}\left|\nabla u_R\(x\)\right|\le C_0
\end{equation}
and 
\begin{equation}\label{Th1Step4Eq5}
\big|R^nf\big(Rx,R^{-\mu}u_R\big)\big|\le\Lambda C_0^{q-p}R^{n-bq-\mu\(q-1\)}\left|x\right|^{-bq-\mu\(q-p\)}u_R^{p-1}
\end{equation}
for a.e. $x\in\R^n\backslash B\(0,1/R\)$. By remarking that 
\begin{equation}\label{Th1Step4Eq6}
n-bq-\mu\(q-1\)=-\frac{\mu\(q-p\)}{p}<0
\end{equation}
and applying the Harnack inequality (see Serrin~\cite{Ser1}), it then follows from \eqref{Th1Step4Eq3}--\eqref{Th1Step4Eq6} that 
\begin{equation}\label{Th1Step4Eq7}
\sup_{B\(0,4\)\backslash B\(0,1/4\)} u_R\le C\hspace{-3pt}\inf_{B\(0,4\)\backslash B\(0,1/4\)} u_R
\end{equation}
for some constant $C>0$ independent of $R>1$. Now we assume by contradiction that \eqref{Th1Step4Eq2} does not hold true. Since $u$ is positive and continuous in $\R^n\backslash\left\{0\right\}$, it then follows from \eqref{Th1Step4Eq7} that there exists a sequence $\(R_k\)_{k\in\N}$ such that
\begin{equation}\label{Th1Step4Eq8}
R_k\to\infty\quad\text{and}\quad\sup_{B\(0,4\)\backslash B\(0,1/4\)} u_{R_k}\to0\quad\text{as }k\to\infty.
\end{equation}
Since $u_{R_k}$ satisfies \eqref{Th1Step4Eq3}--\eqref{Th1Step4Eq6}, by applying the H\"{o}lder estimates of DiBenedetto~\cite{DiB} and Tolksdorf~\cite{Tol}, we obtain that $\(u_{R_k}\)_{k\in\N}$ is bounded in $C^{1,\theta}\(B\(0,3\)\backslash B\(0,1/3\)\)$ for some $\theta\in\(0,1\)$. It then follows from \eqref{Th1Step4Eq8} that up to a subsequence $u_{R_k}\to0$ in $C^1\(B\(0,2\)\backslash B\(0,1/2\)\)$. Let $\eta\in C^1\(\R^n\)$ be a cutoff function such that $\eta\equiv1$ in $B\(0,1/2\)$, $\eta\equiv0$ in $\R^n\backslash B\(0,2\)$ and $0\le\eta\le1$ in $B\(0,2\)\backslash B\(0,1/2\)$. By testing \eqref{Main Problem} with $\eta_k\(x\):=\eta\(x/R_k\)$ and using H\"older's inequality, we obtain
\begin{align}\label{Th1Step4Eq9}
&\int_{\R^n}f\(x,u\)\eta_k\,dx=\int_{\R^n}\left|x\right|^{-ap}\left|\nabla u\right|^{p-2}\left<\nabla u,\nabla\eta_k\right>dx\nonumber\\
&\qquad=\int_{\R^n}\left|x\right|^{-ap}\left|\nabla u_{R_k}\right|^{p-2}\left<\nabla u_{R_k},\nabla\eta\right>dx\nonumber\\
&\qquad\le2^{\left|a\right|p}\left\|\nabla u_{R_k}\right\|^{p-1}_{L^p\(B\(0,2\)\backslash B\(0,1/2\)\)}
\left\|\nabla\eta\right\|_{L^p\(B\(0,2\)\backslash B\(0,1/2\)\)}.
\end{align}
Since up to a subsequence $u_{R_k}\to0$ in $C^1\(B\(0,2\)\backslash B\(0,1/2\)\)$, it follows from \eqref{Th1Step4Eq9} that 
\begin{equation}\label{Th1Step4Eq10}
\int_{\R^n}f\(x,u\)dx=\lim_{k\to\infty}\int_{\R^n}f\(x,u\)\eta_k\,dx=0.
\end{equation}
Since $f\(x,u\)\ge0$ in $\R^n$, it follows from \eqref{Th1Step4Eq10} that $f\(x,u\)\equiv0$ in $\R^n$. Since $u$ is bounded in $\R^n$, by applying a weighted version of Liouville's theorem (see Heinonen, Kilpela\"{i}nen and Martio~\cite{HeiKilMar}*{Theorem~6.10}), we then obtain that $u$ is constant, which is in contradiction with \eqref{Th1Eq3} and the fact that $u$ is positive in $\R^n$. This ends the proof of Step~\ref{Th1Step4}.
\endproof

The next step is as follows:

\begin{step}\label{Th1Step5}
Let $n$, $p$, $a$, $b$, $q$, $\mu$ and $f$ be as in Theorem~\ref{Th1}. Let $u$ be a solution of \eqref{Main Problem} such that $u>0$ and $f\(x,u\)\ge0$ in $\R^n$. Then there exists a constant $\alpha>0$ such that
\begin{equation}\label{Th1Step5Eq1}
\lim_{R\to\infty}\Gamma_R\(u\)=\alpha,\quad\text{where }\Gamma_R\(u\):=\min_{x\in\S^n}\(R^{\mu}u\(Rx\)\).
\end{equation}
\end{step}

\proof
By applying Step~\ref{Th1Step4}, we obtain that
$$\alpha:=\liminf_{R\to\infty}\Gamma_R\(u\)>0.$$
Assume by contradiction that \eqref{Th1Step5Eq1} is not true, namely that
$$\limsup_{R\to\infty}\Gamma_R\(u\)>\alpha.$$
We then obtain that there exist $R_1$, $R_2>0$ such that $R_1<R_2$ and
$$\beta:=\min_{R_1<R<R_2}\(\Gamma_R\(u\)\)<\min\(\Gamma_{R_1}\(u\),\Gamma_{R_2}\(u\)\).$$
It follows that
\begin{equation}\label{Th1Step5Eq2}
\min_{x\in A}\(u\(x\)-w_\beta\(x\)\)=0<\min_{x\in \partial A}\(u\(x\)-w_\beta\(x\)\),
\end{equation}
where $A:=B\(0,R_2\)\backslash B\(0,R_1\)$ and 
$$w_\beta\(x\):=\beta\left|x\right|^{-\mu}\quad\forall x\in\R^n\backslash\left\{0\right\}.$$
By observing that 
$$\divergence\(\left|x\right|^{-ap}\left|\nabla u\right|^{p-2}\nabla u\)\le0=\divergence\(\left|x\right|^{-ap}\left|\nabla w_\beta\right|^{p-2}\nabla w_\beta\)\,\,\,\,\,\text{in }\R^n\backslash\left\{0\right\},$$
we then obtain that \eqref{Th1Step5Eq2} contradicts the strict comparison principle of Serrin~\cite{Ser2}*{Theorem~1}. This ends the proof of Step~\ref{Th1Step5}.
\endproof

We can now end the proof of \eqref{Firstorder} by using Steps~\ref{Th1Step4} and~\ref{Th1Step5}.

\proof[Proof of \eqref{Firstorder}]
Let $u$ be a solution of \eqref{Main Problem} such that $u>0$ and $f\(x,u\)\ge0$ in $\R^n$. By applying Step~\ref{Th1Step5}, we obtain that \eqref{Th1Step5Eq1} holds true. We will prove that \eqref{Firstorder} holds true with $\alpha$ given by \eqref{Th1Step5Eq1}. Assume by contradiction that \eqref{Firstorder} is not true. Then there exists a sequence $\(R_k\)_{k\in\N}$ such that $R_k\to\infty$ and
\begin{equation}\label{Th1Step6Eq1}
\limsup_{k\to\infty}\sup_{x\in\S^n}\left(\left|u_{R_k}\(x\)-\alpha\right|+\left|\nabla u_{R_k}\(x\)+\alpha\mu x\right|\right)>0,
\end{equation}
where $u_{R_k}$ is as in \eqref{Th1Step3Eq2}. We recall that $u_{R_k}$ satisfies \eqref{Th1Step4Eq3}--\eqref{Th1Step4Eq6}. By applying the H\"{o}lder estimates of DiBenedetto~\cite{DiB} and Tolksdorf~\cite{Tol}, we then obtain that for every compact set $\Omega\subset\R^n\backslash\left\{0\right\}$, there exists $\theta_\Omega\in\(0,1\)$ such that $\(u_{R_k}\)_{k\in\N}$ is bounded in $C^{1,\theta_\Omega}\(\Omega\)$ and so there exists a subsequence of $\(u_{R_k}\)_{k\in\N}$ which converges in $C^1\(\Omega\)$ to some function $u_\infty\in C^1\(\R^n\backslash\left\{0\right\}\)$. By passing to the limit as $k\to\infty$ into \eqref{Th1Step4Eq3} and using \eqref{Th1Step4Eq4}--\eqref{Th1Step4Eq6}, we obtain that $u_\infty$ satisfies the equation
$$\divergence\(\left|x\right|^{-ap}\left|\nabla u_\infty\right|^{p-2}\nabla u_\infty\)=0\quad\text{in }\R^n\backslash\left\{0\right\}.$$
By using \eqref{Th1Step5Eq1} and observing that $\Gamma_r\(u_{R_k}\)=\Gamma_{R_kr}\(u\)$, we obtain 
$$\Gamma_r\(u_\infty\)=\lim_{k\to\infty}\Gamma_r\(u_{R_k}\)=\alpha\quad\forall r>0.$$
By another application of the strict comparison principle~\cite{Ser2}*{Theorem~1}, we then obtain 
$$u_\infty\(x\)=\alpha\left|x\right|^{-\mu}\quad\forall x\in\R^n\backslash\left\{0\right\},$$
which contradicts \eqref{Th1Step6Eq1}. This ends the proof of \eqref{Firstorder}.
\endproof

\section{The Kelvin-type transformation}\label{SecKelvin}

This section is devoted to the proof of Proposition~\ref{Pr2} and its application to the proof of \eqref{Th1Eq4}.

\proof[Proof of Proposition~\ref{Pr2}] 
Let $n$, $p$, $a$, $b$, $q$, $\mu$ and $f$ be as in Theorem~\ref{Th1} and $u$ be a weak solution of \eqref{Pr2Eq1}. Assume that there exists a constant $\alpha>0$ such that \eqref{Firstorder} holds true. It follows from \eqref{Firstorder} that if $R$ is chosen large enough, then $u>0$ and $\left|\nabla u\right|>0$ in $\R^n\setminus B\(0,R\)$. Let $r:=1/R$ and $u_\ast$ be the function defined as in \eqref{Pr2Eq2}. Since $u\in C^1\(\R^n \setminus B\(0,R\)\)$, we have $u_\ast\in C^1\big(\overline{B\(0,r\)}\backslash\left\{0\right\}\big)$. In what follows, we will use the notations $\nabla_x$, $\divergence_x$ and $\nabla_y$, $\divergence_y$ for the gradient and divergence with respect to $x$ and $y$, respectively. It follows from \eqref{Firstorder} that $u_\ast$ can be extended to a continuous function in $B\(0,r\)$ and
\begin{align}\label{Pr2Eq6}
\left|\nabla_y u_\ast\right|&=\left|-\mu\left|y\right|^{-\mu-2}uy+\left|y\right|^{-\mu-2}\nabla_x u-2\left|y\right|^{-\mu-4}\<\nabla_x u,y\>y\right|\nonumber\\
&=\left|-\mu\left|x\right|^{\mu}ux+\left|x\right|^{\mu+2}\nabla_x u-2\left|x\right|^{\mu}\<\nabla_x u,x\>x\right|\nonumber\\
&=\smallo\(\left|x\right|\)=\smallo\(\left|y\right|^{-1}\)
\end{align}
as $\left|y\right|\to0$, where $x:=\left|y\right|^{-2}y$. By letting $\gamma$ be as in \eqref{Pr2Eq5} and remarking that $\gamma-2>-n$, it follows from \eqref{Pr2Eq6} that $u_\ast\in H^{1,2}\(B\(0,r\),\left|y\right|^{\gamma}\)$. Furthermore, for every $\varphi\in C^\infty_c\(B\(0,r\)\setminus\left\{0\right\}\)$ straightforward computations give
\begin{align}\label{Pr2Eq7}
&\int_{B\(0,r\)}\left|x\right|^{-ap}\left|\nabla_x u\right|^{p-2}\big<\nabla_x u,\nabla_x\big(\left|x\right|^{-\mu}\varphi\big)\big>dx=\sum_{i,j=1}^n\int_{B\(0,r\)}\left|x\right|^{2n-2-ap}\nonumber\\
&\qquad\times\left|\nabla_x u\right|^{p-2}\partial_{x_i} u\(\delta_{ij} - 2\left|x\right|^{-2} x_i x_j\)\partial_{y_j}\big(\left|x\right|^{-\mu}\varphi\big)dy\nonumber\\
&\quad=\int_{B\(0,r\)}\left|x\right|^{2n-2-ap-\mu}\left|\nabla_x u\right|^{p-2}\<\nabla_x u-2\left|x\right|^{-2}\langle\nabla_xu,x\rangle x,\nabla_y\varphi\>dy\nonumber\\
&\qquad-\mu\int_{B\(0,r\)}\left|x\right|^{2n-2-ap-\mu}\left|\nabla_x u\right|^{p-2}\<\nabla_x u,x\>\varphi\,dy.
\end{align}
We now compute 
\begin{align}\label{Pr2Eq8}
\nabla_x u &= - \mu \left|x\right|^{- \mu-2} u_\ast x  + \left|x\right|^{- \mu - 2} \nabla_y u_\ast - 2 \left|x\right|^{-\mu - 4} \langle \nabla_y u_\ast , x \rangle  x\nonumber\\
&  = - \mu \left|y\right|^{\mu} u_\ast y + \left|y\right|^{\mu+2} \nabla_y u_\ast - 2 \left|y\right|^{\mu} \langle \nabla_y u_\ast , y \rangle y  
\end{align}
and 
\begin{equation}\label{Pr2Eq9}
\langle  \nabla_xu,x \rangle = - \mu \left|y\right|^{\mu} u_\ast -  \left|y\right|^{\mu} \langle \nabla_y u_\ast,y\rangle. 
\end{equation}
It follows from \eqref{Pr2Eq8} and \eqref{Pr2Eq9} that
\begin{equation}\label{Pr2Eq10}
\nabla_x u  - 2 \left|x\right|^{-2} \langle  \nabla_xu,x \rangle x=  \mu  \left|y\right|^{\mu} u_\ast y +   \left|y\right|^{\mu+2} \nabla_y u_\ast
\end{equation}
and 
\begin{equation}\label{Pr2Eq11}
\left|\nabla_x u\right|^2 =  \mu^2  \left|y\right|^{2\mu+2} u_\ast^2  +   2  \mu   \left|y\right|^{2 \mu+2} u_\ast  \langle\nabla_y u_\ast,y \rangle+ \left|y\right|^{2 \mu+4}     \left|\nabla_y u_\ast\right|^2.
\end{equation}
By using \eqref{Pr2Eq9}--\eqref{Pr2Eq11} together with the fact that $\mu\(p-1\)+\(a+1\)p=n$, we obtain
\begin{multline}\label{Pr2Eq12}
\left|x\right|^{2n-2-ap-\mu}\left|\nabla_x u\right|^{p-2}\(\nabla_x u  - 2 \left|x\right|^{-2} \langle  \nabla_xu,x \rangle x\)  \\
= \left|y\right|^{\gamma} \(  \mu^2 u_\ast^2  +  2 \mu  u_\ast \langle \nabla_y u_\ast, y \rangle + \left|\nabla_y u_\ast\right|^2 \left|y\right|^2  \)^{\frac{p-2}{2}}\\
\times\(  \mu  \left|y\right|^{-2} u_\ast y + \nabla_y u_\ast\)
\end{multline}
and 
\begin{multline}\label{Pr2Eq13}
\left|x\right|^{2n-2-ap-\mu}\left|\nabla_x u\right|^{p-2}  \<  \nabla_xu,x \>\\
=-\left|y\right|^{\gamma-2} \( \mu^2 u_\ast^2 + 2 \mu  u_\ast \langle \nabla_y u_\ast, y \rangle+ \left|\nabla_y u_\ast\right|^2 \left|y\right|^2  \)^{\frac{p-2}{2}}\\
\times\(  \mu u_\ast+  \langle \nabla_y u_\ast, y \rangle \).
\end{multline}
Let $\widetilde{A}:B\(0,r\)\times\(0,\infty\)\times\R^n\to\R^n$ and $\widetilde{B}:B\(0,r\)\times\(0,\infty\)\times\R^n\to\R$ be the Caratheodory functions defined by
$$\widetilde{A}\(y,z,\xi\) := \left|y\right|^{\gamma}  E\(y,z,\xi\)^{\frac{p-2}{2}}  \(  \mu  \left|y\right|^{-2} z y + \xi \)-  \mu^{p-1}  \left|y\right|^{\gamma-2} z^{p-1} y$$
and
\begin{multline*}
\widetilde{B}\(y,s,\xi\) :=-\mu \left|y\right|^{\gamma-2}  E\(y,z,\xi\)^{\frac{p-2}{2}}   \(  \mu z+  \langle \xi, y \rangle  \)\\
+    \(p-1\)\mu^{p-1}\left|y\right|^{\gamma-2} z^{p-2}   \langle \xi, y \rangle+\mu^p  \left|y\right|^{\gamma-2} z^{p-1}\\
+ \left|y\right|^{\mu-2n} f\(\left|y\right|^{-2}y, \left|y\right|^{\mu} z\)
\end{multline*} 
for all $\(y,z,\xi\)\in B\(0,r\)\times\(0,\infty\)\times\R^n$, where 
$$E\(y,z,\xi\):=\mu^2 z^2  +  2 \mu  z \langle \xi, y \rangle + \left|\xi\right|^2 \left|y\right|^2.$$
We claim that $u_\ast$ solves the equation
$$-\divergence_y \widetilde{A}\(y,u_\ast,\nabla_y u_\ast\) = \widetilde{B}\(y,u_\ast,\nabla_y u_\ast\)\quad\text{ in } B\(0,r\) \setminus \{0\}.$$
Indeed, by using \eqref{Pr2Eq7}, \eqref{Pr2Eq12} and \eqref{Pr2Eq13} together with the fact that $u$ solves \eqref{Main Problem}, we obtain  
$$\divergence_y \widetilde{A}\(y,u_\ast,\nabla_y u_\ast\)+\widetilde{B}\(y,u_\ast,\nabla_y u_\ast\)=H\(y,u_\ast, \nabla_y u_\ast\)\,\,\,\text{in }B\(0,r\) \setminus \{0\},$$
where 
\begin{multline*}
H\(y,u_\ast, \nabla_y u_\ast\)=-\divergence_y \big(  \mu^{p-1}  \left|y\right|^{\gamma-2}u_\ast^{p-1} y \big)\\
+    \(p-1\)\mu^{p-1}\left|y\right|^{\gamma-2} u_\ast^{p-2}   \langle \nabla_y u_\ast, y \rangle +\mu^p  \left|y\right|^{\gamma-2} u_\ast^{p-1}.  
\end{multline*}
Straightforward computations then give $H\(y,u_\ast, \nabla_y u_\ast\)=0$, thus proving our claim. We now let $A:B\(0,r\)\times\R\times\R^n\to\R^n$ and $B:B\(0,r\)\times\R\times\R^n\to\R$ be the Caratheodory functions defined by
\begin{multline*}
A\(y,z,\xi\):=\widetilde{A}\(y,g_0\(z\),g\(y,\xi\)\)\\
+\mu^{p-2}\left|y\right|^{\gamma}g_0\(z\)^{p-2}\(\(p-2\)\left|y\right|^{-2}\<\xi-g\(y,\xi\),y\>y+\xi-g\(y,\xi\)\)
\end{multline*}
and
$$B\(y,z,\xi\):=\widetilde{B}\(y,g_0\(z\),g\(y,\xi\)\)$$
for all $\(y,z,\xi\)\in B\(0,r\)\times\R\times\R^n$, where
$$g_0\(z\):=\left\{\begin{aligned}&\alpha/2&&\text{if }z<\alpha/2\\&z&&\text{if }\alpha/2\le z\le3\alpha/2\\&3\alpha/2&&\text{if }z>3\alpha/2\end{aligned}\right.$$
and
$$g\(y,\xi\):=\min\(1,\frac{\delta}{\left|y\right|\left|\xi\right|}\)\xi.$$
Here $\delta$ is a positive constant that will be chosen later. It follows from \eqref{Firstorder} and \eqref{Pr2Eq6} that if $R$ is chosen large enough i.e. $r$ is chosen small enough, then $g_0\(u_\ast\)=u_\ast$ and $g\(y,\nabla_y u_\ast\)=\nabla_y u_\ast$ in $B\(0,r\) \setminus \{0\}$ and so $u_\ast$ solves the equation
\begin{equation}\label{Pr2Eq14}
-\divergence_y A\(y,u_\ast,\nabla_y u_\ast\) = B\(y,u_\ast,\nabla_y u_\ast\)\quad\text{ in } B\(0,r\) \setminus \{0\}.
\end{equation}
Furthermore, since $\alpha/2\le g_0\(z\)\le3\alpha/2$ and $\left|g\(y,\xi\)\right|\le\delta/\left|y\right|$, straightforward estimates give 
\begin{multline*}
E\(y,g_0\(z\),g\(y,\xi\)\)^{\frac{p-2}{2}}=\mu^{p-2}g_0\(z\)^{p-2}\\
+\(p-2\)\mu^{p-3}g_0\(z\)^{p-3}\<g\(y,\xi\),y\>+\bigO\(\left|y\right|^2\left|g\(y,\xi\)\right|^2\),
\end{multline*}
which then yields
\begin{multline}\label{Pr2Eq15}
A\(y,z,\xi\)=   \mu^{p-2}  \left|y\right|^{\gamma} g_0\(z\)^{p-2} \( \(p-2\) \left|y\right|^{-2}  \langle \xi, y \rangle y+\xi\)\\
+\bigO\big(\left|y\right|^{\gamma+1}\left|g\(y,\xi\)\right|^2 \big)
\end{multline}
and
\begin{equation}\label{Pr2Eq16}
B\(y,z,\xi\) =\bigO\big(\left|y\right|^{\gamma}\left|g\(y,\xi\)\right|^2\big)+ \left|y\right|^{\mu-2n} f\(\left|y\right|^{-2}y, \left|y\right|^{\mu}g_0\(z\)\)
\end{equation}
uniformly in $\(y,z,\xi\)\in B\(0,r\)\times\R\times\R^n$. By using \eqref{Pr2Eq16} together with \eqref{Th1Eq1} and remarking that $\(b+\mu\)q=\mu q/p+n$ and $\left|g\(y,\xi\)\right|\le\left|\xi\right|$, we then obtain
$$B\(y,z,\xi\)=\bigO\big( \left|y\right|^\gamma\left|g\(y,\xi\)\right|^2  + \left|y\right|^{\gamma'} \big)=\bigO\big( \left|y\right|^\gamma\left|\xi\right|^2  + \left|y\right|^{\gamma'} \big),$$
where $\gamma$ and $\gamma'$ are as in \eqref{Pr2Eq5}. On the other hand, by using \eqref{Pr2Eq15}, we obtain
$$\left|A\(y,z,\xi\)\right|=\bigO\(\left|y\right|^\gamma\left|\xi\right|+\left|y\right|^{\gamma+1}\left|g\(y,\xi\)\right|^2\)=\bigO\(\left|y\right|^\gamma\left|\xi\right|\)$$
and
\begin{align*}
A\(y,z,\xi\)\cdot\xi&=\mu^{p-2}  \left|y\right|^\gamma g_0\(z\)^{p-2} \( \(p-2\) \left|y\right|^{-2}  \langle \xi, y \rangle^2+\left|\xi\right|^2\)\\
&\quad+\bigO\(\left|y\right|^{\gamma+1}\left|g\(y,\xi\)\right|^2\left|\xi\right|\)\\
&\ge\mu^{p-2}  \left|y\right|^\gamma g_0\(z\)^{p-2} \(\min\(p-1,1\)+\bigO\(\delta\)\)\left|\xi\right|^2\\
&\ge C^{-1}\left|y\right|^\gamma\left|\xi\right|^2
\end{align*} 
for some constant $C>0$ independent of $\(y,z,\xi\)\in B\(0,r\)\times\R\times\R^n$, provided we choose $\delta$ small enough. This proves that \eqref{Pr2Eq4} holds true. It remains to show that $u_*$ is a weak solution of \eqref{Pr2Eq3}. Let $\eta\in C^1\(\R^n\)$ be a cutoff function such that $\eta\equiv1$ in $B\(0,1/2\)$, $\eta\equiv0$ in $\R^n\backslash B\(0,1\)$ and $0\le\eta\le1$ in $B\(0,1\)\backslash B\(0,1/2\)$. For every $\varepsilon>0$, let $\eta_\varepsilon\in C^1\(\R^n\)$ be the function defined by $\eta_\varepsilon\(y\)=\eta\(y/\varepsilon\)$ for all $y\in\R^n$. For every $\varphi\in C^\infty_0\(B\(0,r\)\)$, by using $\(1-\eta_\varepsilon\)\varphi$ as a test function, we obtain
\begin{multline}\label{Pr2Eq17}
\int_{B\(0,r\)}\(1-\eta_\varepsilon\)\(A\(y,u_\ast,\nabla u_\ast\)\cdot\nabla\varphi-B\(y,u_\ast,\nabla u_\ast\)\varphi\)dy\\
=\int_{B\(0,r\)}\(A\(y,u_\ast,\nabla u_\ast\)\cdot\nabla\eta_\varepsilon\)\varphi\,dy.
\end{multline}
On the other hand, by using \eqref{Pr2Eq4} and \eqref{Pr2Eq6} together with our definition of $\eta_\varepsilon$, we obtain
\begin{align}\label{Pr2Eq18}
\int_{B\(0,r\)}\eta_\varepsilon A\(y,u_\ast,\nabla u_\ast\)\cdot\nabla\varphi\,dy&=\bigO\(\int_{B\(0,\varepsilon\)}\left|y\right|^{\gamma-1}dy\)\nonumber\\
&=\bigO\(\varepsilon^{\gamma+n-1}\)=\smallo\(1\),
\end{align}
\begin{align}\label{Pr2Eq19}
\int_{B\(0,r\)}\eta_\varepsilon B\(y,u_\ast,\nabla u_\ast\)\varphi\,dy\,&=\bigO\(\int_{B\(0,\varepsilon\)}\big(\left|y\right|^{\gamma-2}+\left|y\right|^{\gamma'}\big)dy\)\nonumber\\
&=\bigO\big(\varepsilon^{\gamma+n-2}+\varepsilon^{\gamma'+n}\big)=\smallo\(1\)
\end{align}
and
\begin{align}\label{Pr2Eq20}
\int_{B\(0,r\)}\(A\(y,u_\ast,\nabla u_\ast\)\cdot\nabla\eta_\varepsilon\)\varphi\,dy&=\bigO\(\frac{1}{\varepsilon}\int_{B\(0,\varepsilon\)}\left|y\right|^{\gamma-1}dy\)\nonumber\\
&=\bigO\(\varepsilon^{\gamma+n-2}\)=\smallo\(1\)
\end{align}
as $\varepsilon\to0$. It follows from \eqref{Pr2Eq17}--\eqref{Pr2Eq20} that
$$\int_{B\(0,r\)}\(A\(y,u_\ast,\nabla u_\ast\)\cdot\nabla\varphi-B\(y,u_\ast,\nabla u_\ast\)\varphi\)dy=0.$$
This proves that $u_*$ is a weak solution of \eqref{Pr2Eq3}. This ends the proof of Proposition~\ref{Pr2}. 
\endproof

We can now end the proof of \eqref{Th1Eq4} and therefore of Theorem~\ref{Th1} by putting together \eqref{Firstorder}, Proposition~\ref{Pr2} and a result of Stredulinski~\cite{Str}.

\proof[End of proof of \eqref{Th1Eq4} and of Theorem~\ref{Th1}]
Let $u$ be a solution of \eqref{Main Problem} such that $u>0$ and $f\(x,u\)\ge0$ in $\R^n$. Let $u_\ast$ be as in \eqref{Pr2Eq2} and $\alpha$ be given by \eqref{Firstorder}. By applying Proposition~\ref{Pr2} together with the H\"older continuity result of Stredulinski~\cite{Str}*{Theorem~3.1.15} (combined with \cite{Str}*{Theorem~2.2.56 and Lemma~3.1.7}), we then obtain that $u_\ast$ is H\"older continuous in $B\(0,r\)$ provided there exist constants $C,\sigma>0$ and $s\in\[1,2\)$ such that 
\begin{equation}\label{Hardy-K-C-N}
\int_{B\(y_0,\rho\)}\left|y\right|^{\gamma'}\left|\varphi\right|^s dy \le C\rho^{s+\sigma-2} \int_{B\(y_0,\rho\)}\left|y\right|^\gamma\left|\nabla\varphi\right|^s dy
\end{equation}
for all $\varphi\in C^\infty_0\(B\(y_0,\rho\)\)$, $y_0\in\R^n$ and $\rho>0$ such that $\overline{B\(y_0,\rho\)}\subset B\(0,r\)$, where $\gamma$ and $\gamma'$ are as in \eqref{Pr2Eq6}. By remarking that $\gamma'>\gamma-2$ and applying a weighted version of the Hardy inequality, which corresponds to the limit case of the Caffarelli--Kohn--Nirenberg inequality~\cite{CafKohnNir}, we obtain that \eqref{Hardy-K-C-N} holds true with $s=2-\sigma$ provided we choose $\sigma$ small enough. Therefore, we obtain that $u_\ast$ is H\"older continuous in $B\(0,r\)$. In particular, we obtain that there exist constants $C,\delta>0$ such that 
\begin{equation}\label{Th1Eq6}
\left|u_\ast\(y\)-\alpha\right|=\left|u_\ast\(y\)-u_\ast\(0\)\right|\le C\left|y\right|^{\delta}\quad\forall y\in B\(0,r\).
\end{equation}
By putting together \eqref{Th1Eq6} with the definition of $u_\ast$, we then obtain that \eqref{Th1Eq4} holds true. This ends the proof of Theorem~\ref{Th1}.
\endproof

\end{document}